\documentclass[psamsfonts,reqno]{amsart}
\usepackage{amssymb,eucal,latexsym,graphics}

\newcommand{\pup}[1]{\textup{(}{#1}\textup{)}}

\newcommand{\conc}{{{}^{\frown}}}
\newcommand{\jirr}{join-ir\-re\-duc\-i\-ble}

\newcommand{\mirr}{meet-ir\-re\-duc\-i\-ble}
\newcommand{\jsd}{join-sem\-i\-dis\-trib\-u\-tive}
\newcommand{\jsdy}{join-sem\-i\-dis\-trib\-u\-tiv\-i\-ty}

\newcommand{\contr}{a contradiction}

\newcommand{\hL}{\widehat{L}}

\newcommand{\CC}{\mathcal{C}}
\newcommand{\MM}{\mathcal{M}}
\newcommand{\PP}{\mathfrak{P}}

\newcommand{\ol}[1]{\overline{#1}}
\newcommand{\set}[1]{\left\{#1\right\}}
\newcommand{\setm}[2]{\set{{#1}\mid{#2}}}
\newcommand{\famm}[2]{\left({#1}\mid{#2}\right)}
\newcommand{\seq}[1]{\langle{#1}\rangle}
\newcommand{\SDn}{\ensuremath{(\mathrm{SD}_{\vee}^n)}}
\newcommand{\SDI}{\ensuremath{(\mathrm{SD}_{\vee}^{\omega})}}
\newcommand{\SP}{\mathbf{S_p}}
\newcommand{\BS}{\mathbb{S}}
\newcommand{\Cls}{\operatorname{Cl}_{\BS}}
\newcommand{\Fls}{\operatorname{Fl}_{\BS}}
\newcommand{\Fis}{\operatorname{Fil}^{\BS}}
\newcommand{\JJd}{\bigvee\nolimits^*}
\newcommand{\jjd}{\mathop{\vee}\nolimits^*}

\newcommand{\Co}{\mathbf{Co}}

\newcommand{\tr}{\vartriangleleft}

\newcommand{\utr}{\trianglelefteq}

\newcommand{\fin}[1]{{#1}^{<\omega}}

\newcommand{\dnw}{\mathbin{\downarrow}}
\newcommand{\upw}{\mathbin{\uparrow}}
\newcommand{\loz}{\mathbin{\lozenge}}

\newcommand{\D}{\mathbin{D}}
\newcommand{\OO}{\mathsf{O}}

\newcommand{\Kc}{K^{\circ}}
\newcommand{\Lc}{L^{\circ}}
\DeclareMathOperator{\Id}{Id}
\DeclareMathOperator{\Fil}{Fil}

\DeclareMathOperator{\J}{J}
\newcommand{\es}{\varnothing}
\newcommand{\FL}{\mathrm{F}_{\mathbf{L}}}

\newcommand{\MI}[2]{{#1}\mathbin{\nearrow}{#2}}
\newcommand{\JI}[2]{{#1}\mathbin{\searrow}{#2}}

\newcommand{\bs}{\boldsymbol{s}}
\newcommand{\bt}{\boldsymbol{t}}

\newcommand{\xs}{\mathbf{s}}
\newcommand{\xt}{\mathbf{t}}

\newcommand{\dv}{\dot{v}}

\numberwithin{equation}{section}

\theoremstyle{plain}

\newtheorem{thm}{Theorem}[section]
\newtheorem{prop}[thm]{Proposition}
\newtheorem{cor}[thm]{Corollary}
\newtheorem{lem}[thm]{Lemma}
\newtheorem{exple}[thm]{Example}
\newtheorem{claim}{Claim}
\newtheorem*{sclaim}{Claim}

\theoremstyle{definition}

\newtheorem{defn}[thm]{Definition}
\newtheorem{pb}{Problem}

\theoremstyle{remark}
\newtheorem{rem}[thm]{Remark}

\newcommand{\qedc}{{\qed}~{\rm Claim~{\theclaim}.}}
\newcommand{\qedsc}{{\qed}~{\rm Claim.}}

\newenvironment{cproof}
{\begin{proof}[Proof of Claim.]}
{\qedc\renewcommand{\qed}{}\end{proof}}

\newenvironment{scproof}
{\begin{proof}[Proof of Claim.]}
{\qedsc\renewcommand{\qed}{}\end{proof}}

\begin{document}
\dedicatory{Dedicated to the memory of Ivan Rival}

\author[F. Wehrung]{Friedrich Wehrung}
  \address{LMNO\\ CNRS UMR 6139\\
           Universit\'e de Caen, Campus II\\
           D\'epartement de Math\'ematiques\\
           B.P. 5186\\
           14032 CAEN Cedex\\
           FRANCE}
  \email{wehrung@math.unicaen.fr}
  \urladdr{http://www.math.unicaen.fr/\~{}wehrung}
\date{\today}

\title[Sublattices of complete lattices]{Sublattices of complete lattices\\
with continuity conditions}

\begin{abstract}
Various embedding problems of lattices into complete lattices are solved. We
prove that for any join-semilattice $S$ with the minimal join-cover refinement
property, the ideal lattice $\Id S$ of $S$ is both algebraic and
dually algebraic. Furthermore, if there are no infinite
$\D$-sequences in $\J(S)$, then $\Id S$ can be embedded into a direct product
of finite lower bounded lattices. We also find a system of infinitary
identities that characterize sublattices of complete, lower continuous, and
\jsd\ lattices. These conditions are satisfied by any (not necessarily
finitely generated) lower bounded lattice and by any locally finite, \jsd\
lattice. Furthermore, they imply M.~Ern\'e's dual staircase distributivity.

On the other hand, we prove that the subspace lattice of any
infinite-dimensional vector space cannot be embedded into any
$\aleph_0$-complete, $\aleph_0$-upper continuous, and $\aleph_0$-lower
continuous lattice. A similar result holds for the lattice of all order-convex
subsets of any infinite chain.
\end{abstract}

\maketitle

\section{Introduction}\label{S:Intro}

It is a classical result that the ideal lattice $\Id L$ of a lattice $L$ is
an algebraic lattice, furthermore, it contains an isomorphic copy of $L$ and
it satisfies the same identities as $L$, see \cite[Lemma~I.4.8]{Grat98}. A
much harder result is that every \emph{modular} lattice embeds, within its
variety, into an algebraic spatial lattice (see Section~\ref{S:Basic} for
precise definitions), see~\cite{HPR}. Say that a lattice is
\emph{bi-algebraic}, if it is both algebraic and dually algebraic. While
investigating lattices of convex subsets, the authors of
\cite{SeWe3,WeSe1} came across the following problem, which is stated as
Problem~5 in \cite{SeWe3}.

\begin{quote}\em
Can every lattice be embedded into some bi-algebraic lattice?
\end{quote}

After having asked several experts in lattice theory, we finally came to the
surprising conclusion that the answer to that question was unknown.
In the present paper, we solve this problem in the negative, see
Section~\ref{S:NonEmb}. More specifically, we prove that both the
lattice of all
subspaces of any infinite-dimensional vector space and the lattice of all
order-convex subsets of any infinite chain cannot be embedded into any
bi-algebraic lattice, see Corollaries~\ref{C:InfVect} and
\ref{C:Co(omega)}.

Nevertheless, it turns out that one can prove many positive
results in this topic that seem to have been unknown until now. We
introduce a new class of lattices, the so-called \emph{fermentable} lattices,
see Definition~\ref{D:Ferm}. The class of fermentable lattices includes the
class of all ideal lattices of finitely generated lower bounded lattices. We
give in Theorem~\ref{T:FermSemil} an alternative proof of a result also
established, with a different method, by M.\,V. Semenova in~\cite{O(P)}:
\emph{Every fermentable lattice can
be embedded into a direct product of finite lower bounded lattices}. This
extends \cite[Corollary~2.2]{AdGo01}, that states that every finitely presented
lower bounded lattice embeds into a direct product of finite lower bounded
lattices. It also extends the result, established
in~\cite{AdDz94} by using \cite{FrNa73}, that the ideal lattice of any
free lattice embeds into some direct product of finite lower bounded
lattices---observe that any direct product of finite
lattices is bi-algebraic. Furthermore, we obtain other related results, such
as: \emph{The ideal lattice of a join-semilattice with the minimal join-cover
refinement property is bi-algebraic}, see Corollary~\ref{C:RP2LC}. A common
extension of lower continuity and \jsdy, called (after M.~Ern\'e) \emph{dual
$*$-distributivity}, is proved for fermentable lattices, see
Corollary~\ref{C:WF2SD+}.

However, these results do not extend to the class of all (not necessarily
finitely generated) lower bounded lattices, see Section~\ref{S:Basic}.
Nevertheless, for those we still obtain partial results, such as the following.

\begin{itemize}
\item[---]
We find a system of infinitary identities characterizing sublattices of
complete, lower continuous, \jsd\ lattices (see Theorem~\ref{T:SubCLCSD+}).
We observe that these `identities' are satisfied by any lower bounded lattice
(see Corollary~\ref{C:LBandSD+lf}) and by any locally finite, \jsd\ lattice
(see Corollary~\ref{C:LocFinSD+}). Furthermore, they imply M. Ern\'e's ``dual
staircase distributivity'' (see Corollary~\ref{C:SDI2DSD}).

\item[---] A finitely generated lower bounded lattice may not be
embeddable into any complete, lower continuous, lower bounded lattice (see
Example~\ref{Ex:NonCplLB}).

\item[---] There exists a locally finite, lower bounded lattice that cannot be
embedded into any complete, upper continuous, \jsd\ lattice (see
Example~\ref{Ex:NonEmbLB}).

\item[---] A lattice has a complete embedding into some complete, lower
continuous, \jsd\ lattice if{f} it satisfies M. Ern\'e's ``dual
$*$-distributivity'' (see Theorem~\ref{T:CpSubCLCSD+}).
\end{itemize}

Some of our results are easy extensions of known results, such as the lower
continuity result proved in Lemma~\ref{L:RP2LC} or the dual $*$-distributivity
result of Corollary~\ref{C:WF2SD+}---still they do not seem to follow right
away from the already existing literature. Some other results of the present
paper seem to be completely new, such as our characterization result of
sublattices of complete, lower continuous, \jsd\ lattices (see
Theorem~\ref{T:SubCLCSD+}). Some patterns of our proof that certain lattices
cannot be embedded into any bi-algebraic lattice (see
Section~\ref{S:NonEmb}) can be found in von~Neumann's classical proof that
the perspectivity relation in a continuous geometry is transitive,
see~\cite{Neum60}.
However, continuous geometries are modular lattices while
our negative results can be applied to non-modular lattices such as those in
Corollary~\ref{C:Co(omega)}. Still, as our results cover
both the modular and the \jsd\ case, putting them in perspective in the present
paper seemed to us worth the effort.

\section{Basic notions}\label{S:Basic}

For a set $X$, we denote by $\fin{X}$ the set of all finite sequences of
elements of $X$, and we denote by $\seq{\bs,\bt}\mapsto\bs\conc\bt$ the
concatenation of finite sequences. We let $\PP(X)$ denote the powerset of $X$.

For partially ordered sets $K$ and $L$, a map
$f\colon K\to L$ is \emph{meet-complete}, if $x=\bigwedge_{i\in I}x_i$ in $K$
implies that $f(x)=\bigwedge_{i\in I}f(x_i)$ in $L$, for all
$x\in K$ and every family $(x_i)_{i\in I}$ of elements of $K$;
``join-complete'' is defined dually. We say that~$f$ is \emph{complete}, if it
is both meet-complete and join-complete.

A lattice $L$ is \emph{lower continuous}, if the equality
   \[
   a\vee\bigwedge X=\bigwedge(a\vee X),
   \]
holds, for any $a\in L$ and any downward directed $X\subseteq L$ such that
$\bigwedge X$ exists. Of course, we put $a\vee X=\setm{a\vee x}{x\in X}$.
Restricting the cardinality of $X$ to be at most $\kappa$, for $\kappa$
either an infinite cardinal or $\infty$, yields \emph{$\kappa$-lower
continuity}. Upper continuity is defined dually.

A lattice is \emph{\jsd}, if it satisfies the quasi-identity
  \[
  x\vee y=x\vee z\Longrightarrow x\vee y=x\vee(y\wedge z).
  \]
For lattices $K$ and $L$, a homomorphism $f\colon K\to L$ is \emph{lower
bounded}, if the preimage under $f$ of any principal filter of $L$ is either
empty or has a least element.
As in \cite{AdGo01}, a lattice~$L$ is \emph{lower
bounded}, if every lattice homomorphism from a finitely generated free lattice
to $L$ is lower bounded. Equivalently, every finitely generated sublattice of
$L$ is lower bounded in the sense of~\cite{FJN}. It is well-known that every
lower bounded lattice is \jsd, see \cite[Theorem~2.20]{FJN}.

For a join-semilattice $S$, we put $S^-=S\setminus\set{0}$, if $S$ has a zero
element, and $S^-=S$, otherwise. We denote by $\J(S)$ the set of
\jirr\ elements of $S$. We say that a subset $\Sigma$ of $S$
is \emph{join-generates} $S$ (resp., \emph{finitely join-generates})~$S$, if
every element of $S$ is a join (resp., a finite join) of elements of
$\Sigma$.

An element $a\in S$ is \emph{compact}, if for every
upward directed subset $X$ of $S$, if $\bigvee X$ is defined and
$a\leq\bigvee X$, then $a\in\dnw X$. We say that $S$ is \emph{algebraic}, if
it is complete and the set of compact elements of $S$ join-generates $S$.
Note that there are other works, such as \cite{AdGS01}, where completeness is
not included in the definition of an algebraic lattice.

We say that~$S$ is \emph{spatial}, if the set
of all completely \jirr\ elements of~$S$ join-generates~$S$.
It is well known that every dually algebraic
lattice is lower continuous---see \cite[Lemma~2.3]{CrDi}
or \cite[Section~1.4]{Comp}, and spatial---see
\cite[Theorem~I.4.22]{Comp} or \cite[Lemma~1.3.2]{Gorb}.

For any $X\subseteq S$, we put
   \begin{align*}
   \dnw X&=\setm{y\in S}{\exists x\in X\text{ such that }y\leq x}\\
   \upw X&=\setm{y\in S}{\exists x\in X\text{ such that }x\leq y}.
   \end{align*}
We abuse notation slightly by putting $\dnw x=\dnw\set{x}$ and
$\upw x=\upw\set{x}$, for all $x\in S$.
For subsets $X$ and $Y$ of $S$, we say that $X$ \emph{refines} $Y$, in notation
$X\ll Y$, if $X\subseteq\dnw Y$.

For $\kappa$ being either a cardinal number or $\infty$, let the prefix
``$\kappa$-'' mean restriction to families of cardinality at most $\kappa$, for
example, a lattice is $\aleph_0$-meet-complete, if every countable subset has a
meet, while it is $\infty$-meet-complete, if it is meet-complete.

\section{Relativizations of the minimal join-cover refinement property}
\label{S:MCRP}

\begin{defn}\label{D:MCRP}
Let $S$ be a join-semilattice and let $\Sigma\subseteq S$.
For an element $a$ of $S^-$, we put
\begin{itemize}
\item
$\CC(a)=\setm{X\subseteq S}
{X\text{ is finite},\ a\notin\dnw X,\text{ and }a\leq\bigvee X}$. The elements
of $\CC(a)$ are called the \emph{nontrivial join-covers of $a$}.

\item
$\MM(a)=\setm{E\in\CC(a)}
{\forall X\in\CC(a),\ X\ll E\text{ implies that }E\subseteq X}$.
The elements of $\MM(a)$ are called the \emph{minimal nontrivial join-covers
of~$a$}, and we put $\MM_\Sigma(a)=\MM(a)\cap\PP(\Sigma)$.
\end{itemize}
Furthermore, we introduce the following properties of the pair
$\seq{S,\Sigma}$:
\begin{itemize}
\item We say that $S$ has the \emph{$\Sigma$-weak minimal join-cover refinement
property}, in short the \emph{$\Sigma$-WMCRP}, if every element of
$\CC(p)$ can be refined by an element of $\MM_\Sigma(p)$, for all
$p\in\Sigma$.

\item We say that $S$ has the \emph{$\Sigma$-minimal join-cover refinement
property}, in short the \emph{$\Sigma$-MCRP}, if it has the $\Sigma$-WMCRP and
$\MM_\Sigma(p)$ is finite, for all $p\in\Sigma$.
\end{itemize}
\end{defn}

Of course, the $\Sigma$-MCRP implies the $\Sigma$-WMCRP. Observe that for
$a\in S$, every element of $\MM_\Sigma(a)$ is an antichain of
$\Sigma\cap\J(S)$. The classical \emph{minimal join-cover refinement property},
in short MCRP, see \cite{FJN}, is the $S$-MCRP. Observe that it implies that
$\J(S)$ finitely join-generates $S$.

The \emph{join-dependency relation} $\D$ on a join-semilattice $S$ is
defined on $\J(S)$ as usual, that is, for $a$, $b\in\J(S)$, the
relation $a\D b$
holds if $a\neq b$ and there exists $c\in S$ such that $a\leq b\vee c$ but
$a\nleq x\vee c$ for all $x<b$. Another useful equivalent definition
is provided
by the following lemma.

\begin{lem}\label{L:EquivD}
Let $S$ be a join-semilattice and let $\Sigma\subseteq\J(S)$. We suppose that
$S$ satisfies the $\Sigma$-WMCRP. Then for all $a$, $b\in\Sigma$, the relation
$a\D b$ holds if{f} there exists $E\in\MM_\Sigma(a)$ such that $b\in E$.
\end{lem}

\begin{proof}
The proof is virtually the same as the one of \cite[Lemma~2.31]{FJN}.
\end{proof}

\section{More lattices with the MCRP}\label{S:moreMCRP}

It is well-known that every free lattice has the MCRP, see
\cite[Chapter~II]{FJN}. Furthermore, every finitely presented lattice has the
MCRP as well, see \cite{Free89}. In this section, we present a few easy common
extensions of these results.

\begin{defn}\label{D:FLP}
For a poset $P$, we denote by $\FL(P)$ the \emph{free lattice on $P$}. This
means that $\FL(P)$ is generated, as a lattice, by (an isomorphic copy of) $P$,
and any order-preserving map from $P$ to any lattice $L$ can be extended to a
unique lattice homomorphism from $\FL(P)$ to $L$.
\end{defn}

\begin{prop}\label{P:FR2mjc}
Every lattice of the form $\FL(P)/{\theta}$, for a poset $P$ and a finitely
generated congruence $\theta$ of $\FL(P)$, has the MCRP.
\end{prop}

\begin{proof}
Put $L=\FL(P)/{\theta}$. We identify $\FL(Q)$ with its canonical image in
$\FL(P)$, for any subposet $Q$ of $P$.

Since $\theta$ is finitely generated,
there are $m<\omega$ and $a_i$, $b_i$, for $i<m$, in $\FL(P)$, such that
$\theta=\bigvee_{i<m}\Theta_{\FL(P)}(a_i,b_i)$. Let $Q_0$ be a
finite subset of $P$ such that $a_i$, $b_i\in\FL(Q_0)$, for all $i<m$.

Let $v\in L$. There exists $\dv\in\FL(P)$ such that $v=[\dv]_{\theta}$. Let
$Q\subseteq P$ be finite such that $Q_0\subseteq Q$ and $\dv\in\FL(Q)$. We put
$\psi=\bigvee_{i<m}\Theta_{\FL(Q)}(a_i,b_i)$, a finitely generated
congruence of $\FL(Q)$. Put $K=\FL(Q)/{\psi}$. There exists a unique lattice
homomorphism $f\colon K\to L$ such that $f([x]_{\psi})=[x]_{\theta}$, for all
$x\in\FL(Q)$.

Put $\Kc=K\cup\set{\OO}$ and $\Lc=L\cup\set{\OO}$, for a new zero
element $\OO$. Extend $f$ to a homomorphism from $\Kc$ to $\Lc$ by putting
$f(\OO)=\OO$. Since $Q$ is finite, it is possible to define a map
$g_0\colon P\to\Kc$ by the rule
  \[
  g_0(p)=\bigvee\famm{[q]_\psi}{q\in Q,\ q\leq p},\text{ for all }p\in P,
  \]
with the convention that $\bigvee\es=\OO$. Since $g_0$ is order-preserving, it
extends to a unique lattice homomorphism $g\colon\FL(P)\to\Kc$. Observe that
$g(x)=[x]_\psi$, for all $x\in\FL(Q)$. In particular, $g(a_i)=g(b_i)$, for all
$i<m$, thus there exists a unique lattice homomorphism $h\colon\Lc\to\Kc$ such
that $h(\OO)=\OO$ and $h([x]_\theta)=g(x)$, for all $x\in\FL(P)$.

\begin{sclaim}
The following assertions hold.
\begin{enumerate}
\item $hf(y)=y$, for all $y\in\Kc$; thus $f$ is an embedding.

\item The inequality $fh(y)\leq y$ holds, for all $y\in\Lc$.
\end{enumerate}
\end{sclaim}

\begin{scproof}
(i) In the nontrivial case $y\neq\OO$, we can write $y=[t]_\psi$, for
some $t\in\FL(Q)$. Then $hf(y)=h([t]_\theta)=g(t)=[t]_\psi=y$.

(ii) For all $p\in P$, we compute:
   \[
   fh([p]_\theta)=fg(p)=\bigvee\famm{f([q]_\psi)}{q\leq p\text{ in }Q}
   =\bigvee\famm{[q]_\theta}{q\leq p\text{ in }Q}\leq[p]_\theta.
   \]
Since $\setm{[p]_\theta}{p\in P}$ generates $\FL(P)$, the conclusion follows.
\end{scproof}

Set $u=[\dv]_{\psi}$. Hence $v=f(u)$ and, by (i) of the Claim above, $u=h(v)$.
Since~$K$ is a finitely presented lattice, it has the MCRP, see \cite{Free89}.
Let $I_l$, for $l<n$, denote the minimal nontrivial join-covers of $u$ in
$K$. Since $v=f(u)$, the set $J_l=f[I_l]$ is a nontrivial join-cover
of $v$, for
all $l<n$. Hence, to conclude the proof, it suffices to establish that every
nontrivial join-cover $J$ of $v$ in $L$ is refined by some $J_l$. Observe first
that $u=h(v)\leq\bigvee h[J]$. If $u\leq h(t)$, for some $t\in J$,
then, by (ii)
of the Claim above, $v=f(u)\leq fh(t)\leq t$, a contradiction; hence,
$h[J]$ is a nontrivial join-cover of $u$ in $K$, thus there exists $l<n$ such
that $I_l\ll h[J]$, whence $J_l=f[I_l]\ll fh[J]$, so, again by (ii) of the
Claim above, $J_l\ll J$, which concludes the proof.
\end{proof}

We observe the following immediate consequence of \cite[Lemma~5.3]{FJN}.

\begin{prop}\label{P:LB2mjc}
Let $f\colon L\twoheadrightarrow K$ be a lower bounded homomorphism of
lattices. If~$L$ has the MCRP, then so does $K$.
\end{prop}

\begin{cor}\label{C:LB2mjc}
Every lower bounded homomorphic image of $\FL(P)/{\theta}$ has the MCRP, for
any poset $P$ and any finitely generated congruence $\theta$ of $\FL(P)$.
\end{cor}

\begin{defn}\label{D:FinDef}
A lattice $L$ is \emph{finitely defined}, if it is defined by finitely many
relations (within the class of all lattices).
\end{defn}

It is clear that finitely defined lattices are exactly the quotients of free
lattices by finitely generated congruences. Hence, Corollary~\ref{C:LB2mjc}
applies to finitely defined lattices.

\section{Leavens and fermentable lattices}\label{S:FermLatt}

An \emph{infinite $\D$-sequence} of a join-semilattice $L$ is a sequence
$(a_n)_{n<\omega}$ of elements of $\J(L)$ such that $a_n\D a_{n+1}$ for all
$n<\omega$.

\begin{defn}\label{D:Ferm}
Let $L$ be a join-semilattice. A subset $\Sigma$ of $\J(L)$ is a \emph{leaven}
of $L$, if the following statements hold:
\begin{enumerate}
\item $\Sigma$ join-generates $L$.

\item $L$ satisfies the $\Sigma$-MCRP.

\item There is no infinite $\D$-sequence of elements of $\Sigma$ in $L$.
\end{enumerate}
We say that $L$ is \emph{fermentable}, if it has a leaven.
\end{defn}

This terminology is inspired from P. Pudl\'ak and J. T\r{u}ma's
beautiful designation as ``finitely fermentable'' (see \cite{PuTu74}) those
lattices that are nowadays called ``finite lower bounded''. In particular,
a finite lattice is fermentable if{f} it is lower bounded.

Every free lattice and every finitely generated lower bounded
lattice is fermentable (see \cite{FJN}). Furthermore,
the ideal lattice of the free lattice $\FL(X)$ on any nonempty set $X$ is,
by Corollary~\ref{C:FermSemil}, fermentable, but it is not lower bounded in
case $X$ has at least three elements (see \cite{AdGo01}).

The following result that every fermentable lattice can be embedded into a
direct product of finite lower bounded lattices is also established by M.\,V.
Semenova in~\cite{O(P)}.

\begin{thm}\label{T:FermSemil}
Every fermentable join-semilattice $L$ has a meet-complete join-embedding into
some direct product of finite lower bounded lattices.
\end{thm}

Observe that every direct product of finite lower bounded lattices is
bi-algebraic and fermentable.

\begin{proof}
Let $\Sigma$ be a leaven of $L$. We denote by $\utr$ the reflexive, transitive
closure of the join-dependency relation on $\Sigma$ (see Lemma~\ref{L:EquivD}),
and we put
  \[
  \Sigma_p=\setm{q\in\Sigma}{p\utr q},\text{ for all }p\in\Sigma.
  \]

\setcounter{claim}{0}
\begin{claim}\label{Cl:SpFin}
The set $\Sigma_p$ is finite, for all $p\in\Sigma$.
\end{claim}

\begin{cproof}
Define $T$ as the set of all finite sequences $\seq{p_0,p_1,\dots,p_n}$ of
elements of $\Sigma$ such that $p_i\D p_{i+1}$, for all $i<n$. Then
$T$, endowed with the initial segment ordering, is a tree. Furthermore, since
$L$ satisfies the $\Sigma$-MCRP, $T$ is finitely branching. Since there is no
infinite $\D$-sequence in $\Sigma$, the tree $T$ has no infinite branch, thus,
by K\"onig's Theorem, every connected component of $T$ is finite. The
conclusion
follows immediately.
\end{cproof}

For any $p\in\Sigma$, denote by $L_p$ the set of all joins of elements
of $\Sigma_p$, with a new zero element $\OO$ added as the join of the empty
set.

\begin{claim}\label{Cl:LpLB}
The lattice $L_p$ is finite lower bounded, for all $p\in\Sigma$.
\end{claim}

\begin{cproof}
The finiteness of $L_p$ follows from Claim~\ref{Cl:SpFin}.
Moreover, $\Sigma_p$ is contained in $\J(L)\cap L_p$, thus in $\J(L_p)$. Since
every element of $L_p$ is a join of elements of $\Sigma_p$, it follows that
$\J(L_p)=\Sigma_p$.

For all $q$, $r\in\Sigma_p$, the relations $q\D_{L_p}r$ and $q\D_{L}r$ are
equivalent. Thus, $L_p$ does not have $\D$-cycles. Since $L_p$ is finite, it is
lower bounded.
\end{cproof}

We put $\Phi_p(x)=\dnw x\cap\Sigma_p$ and $\varphi_p(x)=\bigvee\Phi_p(x)$, for
all $p\in\Sigma$ and $x\in L$. In particular, $\varphi_p$ is a map from $L$ to
$L_p$. It is obvious that $q\leq x$ if{f} $q\leq\varphi_p(x)$, for all
$q\in\Sigma_p$ and all $x\in L$; hence
$\varphi_p$ is meet-complete. Let $a$, $b\in L$ and let $q\in\Sigma_p$ with
$q\leq\varphi_p(a\vee b)$, we prove that
$q\leq\varphi_p(a)\vee\varphi_p(b)$. This is obvious if either $q\leq a$ or
$q\leq b$, so suppose that $q\nleq a,b$. Since $q\leq a\vee b$, there exists
$I\in\MM_\Sigma(q)$ such that $I\ll\set{a,b}$. {}From $I\in\MM_\Sigma(q)$ it
follows that $I\subseteq\Sigma_p$, thus, since $I\ll\set{a,b}$, the relation
$I\subseteq\Phi_p(a)\cup\Phi_p(b)$ holds, whence
$q\leq\bigvee I\leq\varphi_p(a)\vee\varphi_p(b)$. It follows that
$\varphi_p$ is
a join-homomorphism from $L$ to $L_p$.

Hence, the map $\varphi\colon L\to\prod_{p\in\Sigma}L_p$ defined by
the rule $\varphi(x)=(\varphi_p(x))_{p\in\Sigma}$ is a join-homomorphism, and
it is meet-complete. Furthermore, for $a$, $b\in L$ such that
$a\nleq b$, there exists $p\in\Sigma$ such that $p\leq a$ and $p\nleq b$, thus
$p\leq\varphi_p(a)$ while $p\nleq\varphi_p(b)$; whence
$\varphi(a)\nleq\varphi(b)$. Therefore, $\varphi$ is an order-embedding.
\end{proof}

\begin{rem}\label{Rk:DecFerm}
The proof above shows, in fact, that all the lattice homomorphisms
$\varphi_p\colon L\to L_p$, for $p\in\Sigma$, are \emph{lower bounded}.
\end{rem}

Every finitely generated lower bounded lattice $L$ is fermentable, with
leaven $\J(L)$ (see \cite[Theorem~2.38]{FJN}); hence, by
Theorem~\ref{T:FermSemil}, it embeds into a direct product of finite lower
bounded lattices. Therefore, every finitely generated lower bounded lattice
belongs to the quasivariety $\mathbf{Q}(\mathcal{LB}_{\mathrm{fin}})$ generated
by all finite lower bounded lattices. As every quasivariety is closed under
direct limits, this gives another proof of the result, first established in
\cite[Theorem~2.1]{AdGo01}, that every lower bounded lattice belongs to
$\mathbf{Q}(\mathcal{LB}_{\mathrm{fin}})$.

\begin{cor}\label{C:FermSemil}
Let $S$ be a join-semilattice with the MCRP and no infinite $\D$-sequence
of \jirr\ elements. Then the ideal lattice $\Id S$ is fermentable; thus it
embeds into a direct product of finite lower bounded lattices.
\end{cor}

\begin{proof}
It is straightforward to verify that $\Sigma=\setm{\dnw p}{p\in\J(S)}$ is a
leaven of $L=\Id S$; whence $L$ is fermentable. The conclusion follows
immediately from Theorem~\ref{T:FermSemil}.
\end{proof}

In particular, it follows from Corollary~\ref{C:FermSemil} that $\Id L$ is
fermentable, for every lattice $L$ which is either free or finitely generated
lower bounded. It cannot be extended to
arbitrary lower bounded lattices, for any direct product of finite lower
bounded lattices is complete, upper continuous, and \jsd, while, on the other
hand, the (locally finite, lower bounded) lattice of Example~\ref{Ex:NonEmbLB}
cannot be embedded into any complete, upper continuous, \jsd\ lattice.
Furthermore, by Example~\ref{Ex:FilFL(3)}, Corollary~\ref{C:FermSemil}
cannot be extended to the \emph{filter lattice} $\Fil L$ of $L$.

\section{Lower continuity}\label{S:WF2LC}

Our main lemma is the following, very similar in statement and in proof to
\cite[Theorem~2.25]{FJN}.

\begin{lem}\label{L:RP2LC}
Let $L$ be a join-semilattice and let $\Sigma$ be a join-generating subset of
$\J(L)$ such that~$L$ has the $\Sigma$-MCRP. Then $L$ is lower continuous.
\end{lem}

\begin{proof}
Let $a\in L$ and let $X\subseteq L$ be a downward directed subset admitting a
meet, $b=\bigwedge X$. We prove the equality $a\vee b=\bigwedge(a\vee X)$.
Since $\Sigma$ join-generates $S$, it suffices to prove that for any
$p\in\Sigma$ such that $p\leq a\vee x$ for all $x\in X$, the inequality
$p\leq a\vee b$ holds. This is trivial in case either $p\leq a$ or $p\leq x$
for all $x\in X$, so suppose that this does not occur; let $x_0\in X$
such that $p\nleq x_0$, and put $X'=\setm{x\in X}{x\leq x_0}$. We put
$\mu(x)=\setm{E\in\MM_\Sigma(p)}{E\ll\set{a,x}}$, for all $x\in X'$.
Observe that $\mu(x)\neq\es$ (because of the $\Sigma$-MCRP)
and that $x\leq y$ implies that $\mu(x)\subseteq\mu(y)$, for all $x\leq y$ in
$X'$. Since $\MM_\Sigma(p)$ is finite (because of the $\Sigma$-MCRP), the
intersection of all $\mu(x)$, for $x\in X'$, is nonempty;
pick an element $E$ in this set. Since $E\ll\set{a,x}$, for all $x\in X'$, the
relation $E\ll\set{a,b}$ holds, whence $p\leq\bigvee E\leq a\vee b$.
\end{proof}

The following lemma is folklore.

\begin{lem}\label{L:AlgLC}
Every algebraic and lower continuous lattice is dually algebraic.
\end{lem}

\begin{proof}
Let $L$ be an algebraic and lower continuous lattice. Since $L$ is dually
spatial (see \cite[Theorem~I.4.22]{Comp}), it suffices to prove that every
completely \mirr\ element $u$ of $L$ is dually compact. Let $X$ be a downward
directed subset of~$L$ such that $\bigwedge X\leq u$. Suppose that $x\nleq u$
for all $x\in X$. So $u^*$, the unique upper cover of $u$, lies below
$u\vee x$, for all $x\in X$; whence, by the lower continuity of~$L$,
$u^*\leq u\vee\bigwedge X=u$, a contradiction. Hence, $u$ is dually compact.
\end{proof}

\begin{cor}\label{C:RP2LC}
Let $S$ be a join-semilattice satisfying the MCRP. Then the ideal lattice
$\Id S$ of $S$ is bi-algebraic.
\end{cor}

\begin{proof}
Apply Lemmas~\ref{L:RP2LC} and \ref{L:AlgLC} to $L=\Id S$ and
$\Sigma=\setm{\dnw p}{p\in\J(S)}$.
\end{proof}

\begin{cor}\label{C:Id(LB)LC}
Let $P$ be a poset, let $\theta$ be a finitely generated congruence of
$\FL(P)$, and let $L$ be a lower bounded homomorphic image of
$\FL(P)/{\theta}$. Then $\Id L$ is bi-algebraic.
\end{cor}

In particular, we obtain the following corollary.

\begin{cor}\label{C:Id(LBFD)LC}
Let $L$ be a lower bounded homomorphic image of a finitely defined
lattice. Then $\Id L$ is bi-algebraic.
\end{cor}

Corollary~\ref{C:Id(LBFD)LC} cannot be extended to the class of all lower
bounded lattices (see Section~\ref{S:Basic}). In fact,
Corollary~\ref{C:Id(LBFD)LC} does not even extend to \emph{Boolean} lattices.

\begin{prop}\label{L:InfBA}
Let $B$ be an infinite Boolean lattice. Then $\Id B$ is not lower continuous.
\end{prop}

\begin{proof}
Let $(a_n)_{n<\omega}$ be a strictly increasing sequence of elements of $B$.
Let $A$ denote the ideal of $B$ generated by $\setm{a_n}{n<\omega}$ and let
$B_n$ denote the principal ideal generated by $\neg a_n$, for all $n<\omega$.
Observe that the sequence $(B_n)_{n<\omega}$ is (strictly) decreasing. Then
the top element $1$ belongs to $A\vee B_n$, for every $n<\omega$, but it does
not belong
to $A\vee\bigcap_{n<\omega}B_n$.
\end{proof}

\section{Dual $*$-distributivity}\label{S:*distr}

Our $*$ operation is the dual of the one considered in \cite{Erne88,Rein95} and
\cite[Section~5.6]{FJN}.

\begin{defn}\label{L:Def*}
For a lattice $L$, we consider the lattice $L\cup\set{1}$ obtained by adding a
new largest element $1$ to $L$. For $a\in L$, we define inductively an element
$a*\bs$ of $L\cup\set{1}$, for $\bs\in\fin{L}$, as follows:
   \begin{align*}
   a*\es&=1;\\
   a*(\bs\conc\seq{b})&=a\vee(b\wedge(a*\bs)),\text{ for all }\bs\in\fin{L}
   \text{ and }b\in L.
   \end{align*}
For a subset $B$ of $L$, we put
   \[
   a*B=\setm{a*\bs}{\bs\in\fin{B}}.
   \]
We say that $L$ is \emph{dually $*$-distributive}, if whenever $a\in L$
and $B\subseteq L$, if $\bigwedge B$ exists,
then $\bigwedge(a*B)$ exists and
   \[
   \bigwedge(a*B)=a\vee\bigwedge B.
   \]
Restricting the cardinality of $B$ to be at most $\kappa$, for a given cardinal
number $\kappa$, defines \emph{dual $\kappa$-$*$-distributivity} \pup{see
Section~\ref{S:Basic}}. We say that $L$ is
\begin{itemize}
\item[---] \emph{dually staircase distributive}, if it is dually
$n$-$*$-distributive, for every positive integer $n$,

\item[---] \emph{dually zipper distributive}, if it is dually
$2$-$*$-distributive.
\end{itemize}
\end{defn}

The following is essentially due to M. Ern\'e, see the proof of
\cite[Proposition~2.12]{Erne88}.

\begin{prop}\label{P:StD}
Let $L$ be a lattice and let $\kappa$ be either an infinite cardinal number
or~$\infty$. Consider the following statements:
\begin{enumerate}
\item $L$ is dually $\kappa$-$*$-distributive.

\item $L$ is $\kappa$-lower continuous and dually staircase distributive.

\item $L$ is $\kappa$-lower continuous and dually zipper distributive.

\item $L$ is $\kappa$-lower continuous and \jsd.
\end{enumerate}
Then \textup{(i)} and \textup{(ii)} are equivalent, they imply \textup{(iii)},
which implies \textup{(iv)}. Furthermore, if $L$ is
$\kappa$-meet-complete, then all four statements are equivalent.
\end{prop}

It will turn out that (i), (ii), and (iii) are, in fact, \emph{equivalent}, see
Corollary~\ref{C:ZiStD}.

The following fact is obvious.

\begin{lem}\label{L:SubL*}
For a sublattice $K$ of a lattice $L$, the following statements hold.
\begin{enumerate}
\item If $L$ is dually staircase distributive \pup{resp., dually zipper
distributive}, then so is $K$.

\item If the inclusion map from $K$ into $L$ is meet-complete and $L$ is dually
$*$-distributive, then so is $K$.
\end{enumerate}
\end{lem}

The following result extends \cite[Theorem~5.66]{FJN}.

\begin{cor}\label{C:WF2SD+}
Every fermentable lattice is dually $*$-distributive.
\end{cor}

\begin{proof}
It follows from Theorem~\ref{T:FermSemil} that $L$ has a meet-complete
lattice embedding into some direct product $\ol{L}$ of finite lower bounded
lattices. Of course, $\ol{L}$ is complete, lower continuous, and \jsd, hence,
by Proposition~\ref{P:StD}, it is dually $*$-distributive. Therefore, by
Lemma~\ref{L:SubL*}, $L$ is also dually $*$-distributive.
\end{proof}

\begin{cor}\label{C:WF2SD+2}
Let $S$ be a join-semilattice satisfying the following assumptions:
\begin{enumerate}
\item $S$ has the MCRP.

\item There are no infinite $\D$-sequences in $\J(S)$.
\end{enumerate}
Then $\Id S$ is dually $*$-distributive; in particular, it is \jsd.
\end{cor}

\begin{proof}
Apply Corollary~\ref{C:WF2SD+} to $L=\Id S$ and
$\Sigma=\setm{\dnw p}{p\in\J(S)}$.
\end{proof}

It is well-known that free lattices have the MCRP and have no infinite
$\D$-sequences, see \cite[Chapter~II]{FJN}. Hence, we obtain the following
consequence.

\begin{cor}\label{C:Id(LB)}
Let $L$ be a lower bounded homomorphic image of a free lattice. Then $\Id L$ is
dually $*$-distributive; in particular, it is \jsd.
\end{cor}

Compare with Corollary~\ref{C:Id(LBFD)LC}.
Observe that we do not require $L$ to be finitely generated. By
Corollary~\ref{C:FermSemil}, $\Id L$ satisfies many other quasi-identities than
\jsdy, namely, all those quasi-identities that hold in all finite lower bounded
lattices, see Theorem~4.2.8 and Corollary~5.5.8 in \cite{Gorb}.

\section{The axiom $\SDI$}\label{S:AxSDI}

Let $L$ be a lattice. For any $\bs=\seq{a_0,\ldots,a_{n-1}}\in\fin L$, we
put $n=|\bs|$, and, if $n>0$, we put $\bs_*=\seq{a_0,\ldots,a_{n-2}}$ and
$e(\bs)=a_{n-1}$.
Furthermore, we define inductively $\bs\loz x$, for
$\bs\in\fin L$ and $x\in L$:
  \begin{align}
  \es\loz x&=x;\label{Eq:escircx}\\
  (\bs\conc\seq{a})\loz x&=
  \begin{cases}\label{Eq:bscirx}
  a\vee(\bs\loz x),&\text{if }|\bs|\text{ is even},\\
  a\wedge(\bs\loz x),&\text{if }|\bs|\text{ is odd}.
  \end{cases}
  \end{align}
We shall need in Section~\ref{S:SubCLCSD+} the following simple lemma.

\begin{lem}\label{L:circjm}
There are maps $j\colon L\times\fin L\to\fin L$ and
$m\colon L\times L\times\fin L\to\fin L$ such that the following equalities
hold for all $u$, $x$, $y\in L$ and all $\bs\in\fin L$ with~$y\leq x$:
  \begin{align*}
  u\vee(\bs\loz x)&=j(u,\bs)\loz x;\\
  u\wedge(\bs\loz x)&=m(u,y,\bs)\loz x.
  \end{align*}
\end{lem}

\begin{proof}
We define the maps $j$ and $m$ inductively, by
\[
  \begin{aligned}
  j(u,\es)=\seq{u}&\quad\text{and}\quad m(u,y,\es)=\seq{y,u};\\
  j(u,\bs)=\bs\conc\seq{u}&\quad\text{and}\quad 
  m(u,y,\bs)=\bs_*\conc\seq{u\wedge e(\bs)},
  \text{ if }|\bs|\text{ is nonzero even};\\
  j(u,\bs)=\bs_*\conc\seq{u\vee e(\bs)}&\quad\text{and}\quad 
  m(u,y,\bs)=\bs\conc\seq{u},\text{ if }|\bs|\text{ is odd.}
  \end{aligned}
\]
It is straightforward to verify that these maps satisfy the required
conditions.
\end{proof}

\begin{defn}\label{D:SDinfty}
We say that a lattice $L$ satisfies \SDI, if the following equality holds, for
all $\bs\in\fin L$ and all $a$, $b$, $c\in L$:
  \begin{equation}\label{Eq:SDI}
  \bs\loz(a\vee(b\wedge c))=\bigwedge
  \setm{\bs\loz(a*\bt)}{\bt\in\fin{\set{b,c}}}.
  \end{equation}
\end{defn}

It is not hard to verify that \SDI\ is a weakening of the identity \SDn\
considered in \cite[Section~4.2]{JiRo}. A very closely related notion is the
(dual) \emph{$m$-zipper distributivity} considered in~\cite{Erne88}.

\begin{prop}\label{P:SD+2SDI}
Let $L$ be a lattice. Consider the following statements:
\begin{enumerate}
\item $L$ satisfies \SDI.

\item $L$ is dually zipper distributive.
\end{enumerate}
Then \textup{(i)} implies \textup{(ii)}. Furthermore, if $L$ is
$\aleph_0$-lower continuous, then \textup{(i)} and \textup{(ii)} are
equivalent.
\end{prop}

\begin{proof}
It is obvious that dual zipper distributivity of $L$ is equivalent to the
satisfaction of \eqref{Eq:SDI} for all elements $a$, $b$, $c$ of $L$ and for
$\bs=\es$, and thus it follows from \SDI. If $L$ is $\aleph_0$-lower
continuous, then it is easy to establish, by induction on the length of $\bs$,
the equality
  \[
  \bs\loz\bigwedge X=\bigwedge
  \setm{\bs\loz x}{x\in X},
  \]
for all $\bs\in\fin{L}$ and every (at most) countable downward directed subset
$X$ of $L$. Under such conditions, dual zipper distributivity obviously implies
\SDI.
\end{proof}

Part of the conclusion of Proposition~\ref{P:SD+2SDI} will be strengthened in
Corollary~\ref{C:SDI2DSD}.

We leave to the reader the easy proof of the following preservation
result.

\begin{prop}\label{P:PresSDI}
The following statements hold.
\begin{enumerate}
\item Every sublattice of a lattice satisfying \SDI\ satisfies \SDI.

\item Every directed union of a family of lattices satisfying \SDI\ satisfies
\SDI.
\item Every image of a lattice satisfying \SDI\ under a
$\aleph_0$-meet-complete lattice homomorphism satisfies \SDI.

\item Every lower bounded homomorphic image of a lattice satisfying \SDI\
satisfies \SDI.
\end{enumerate}
\end{prop}
(Since every lower bounded surjective homomorphism is
meet-complete, item~(iii) trivially implies item (iv).)

As an immediate consequence of Proposition~\ref{P:PresSDI}(i,ii), we obtain the
following.

\begin{cor}\label{C:FinGenSDI}
A lattice $L$ satisfies \SDI\ if{f} every finitely generated sublattice of~$L$
satisfies \SDI.
\end{cor}

\begin{cor}\label{C:DNoethSDI}
Every directed union of fermentable lattices satisfies \SDI.
\end{cor}

\begin{proof}
By Corollary~\ref{C:WF2SD+}, every fermentable lattice is dually
$*$-distributive, hence, by Propositions~\ref{P:StD} and \ref{P:SD+2SDI}, it
satisfies \SDI. The conclusion follows from Proposition~\ref{P:PresSDI}(ii).
\end{proof}

\section{Lower continuous lattices of filters}
\label{S:ConvFil}

\begin{defn}\label{D:CovNot}
A \emph{notion of convergence} on a lattice $L$ is a set $\BS$ of
subsets of $L$
satisfying the following conditions:
\begin{itemize}
\item[(S1)] Every element $X$ of $\BS$ is downward directed, furthermore,
$\bigwedge X$ exists.

\item[(S2)] For all $X\in\BS$ and all $a\in L$, the subset
$a\vee X=\setm{a\vee x}{x\in X}$ belongs to~$\BS$, and
$\bigwedge(a\vee X)=a\vee\bigwedge X$.

\item[(S3)] For all $X\in\BS$ and all $a\in L$, the subset
$a\wedge X=\setm{a\wedge x}{x\in X}$ belongs to~$\BS$. (\emph{Observe that
necessarily, $\bigwedge(a\wedge X)=a\wedge\bigwedge X$}.)
\end{itemize}

We say that $\BS$ is \emph{special}, if it satisfies the following condition,
that involves the~$*$ operation introduced in Section~\ref{S:*distr}:
\begin{itemize}
\item[(S4)] The subset
$\setm{a*\bt}{\bt\in\fin{\set{b,c}}\setminus\set{\es}}$ belongs to
$\BS$ and its meet is $a\vee(b\wedge c)$, for all $a$, $b$, $c\in L$.
\end{itemize}
\end{defn}

As usual, we say that a \emph{filter} of a lattice $L$ is a (possibly empty)
upper subset of~$L$, closed under finite meets.

\begin{defn}\label{D:BSClos}
Let $\BS$ be a set of subsets of a lattice $L$ satisfying (S1).
We say that a subset $A$ of $L$ is \emph{$\BS$-closed}, if $X\subseteq A$
implies that $\bigwedge X\in A$, for all $X\in\BS$.

For a subset $A$ of $L$, we denote by $\Cls(A)$ the least $\BS$-closed subset
of $L$ containing~$A$, and by $\Fls(A)$ the least $\BS$-closed filter
containing $A$.
\end{defn}

\emph{For the remainder of Section~\textup{\ref{S:ConvFil}}, let $\BS$ be a
notion of convergence on a lattice~$L$.}
\smallskip

For any $a\in L$ and any $X\subseteq L$, we put
  \begin{align}
  \MI{X}{a}&=\setm{y\in L}{a\wedge y\in X},\label{Eq:arghX}\\
  \JI{X}{a}&=\setm{y\in L}{a\vee y\in X}.\label{Eq:Xsda}
  \end{align}
The proof of the following lemma is a straightforward application of (S2)
and~(S3).

\begin{lem}\label{L:avnX}
If a subset $X$ of $L$ is $\BS$-closed, then so are $\MI{X}{a}$ and
$\JI{X}{a}$, for any $a\in L$.
\end{lem}

It is obvious that $\Cls(A)$ is contained in $\Fls(A)$, for any $A\subseteq L$.
The following lemma gives us an important case where the two closures
are equal.

\begin{lem}\label{L:MonClass}
The equality $\Cls(A)=\Fls(A)$ holds, for any filter $A$ of $L$.
\end{lem}

\begin{proof}
Put $B=\Cls(A)$. For any $a\in A$, the subset $\MI{B}{a}$ is, by
Lemma~\ref{L:avnX}, $\BS$-closed, and it contains $A$ (because $A$ is closed
under finite meets); thus $B\subseteq\MI{B}{a}$. Hence
$A\subseteq\MI{B}{b}$, for all $b\in B$. But $\MI{B}{b}$ is
$\BS$-closed, thus $B\subseteq\MI{B}{b}$, for all $b\in B$; that is, $B$
is closed under finite meets.

Furthermore, for any $x\in L$, the subset $\JI{B}{x}$ is, by
Lemma~\ref{L:avnX},
$\BS$-closed, but it contains $A$ (because $A$ is an upper subset of $L$), thus
it contains $B$. Hence $B$ is an upper subset of $L$. Therefore, $B$ is a
filter of $L$, but it is $\BS$-closed, thus it contains $\Fls(A)$; whence
$B=\Fls(A)$.
\end{proof}

We denote by $\Fis L$ the set of all $\BS$-closed filters $X$ of $L$ such that
if $L$ has a unit element, say, $1_L$, then $1_L\in X$. (We take this
precaution in order to ensure that the canonical embedding from $\seq{L,\leq}$
into $\seq{\Fis L,\supseteq}$ preserves the empty meet.) Since $\Fis L$ is a
closure system in the powerset of~$L$, the poset $\seq{\Fis L,\subseteq}$ is a
complete lattice. We shall order $\Fis L$ by \emph{reverse} inclusion. The
proof of the following lemma is obvious.

\begin{lem}\label{L:CanEmbS}
The map $x\mapsto\upw x$ defines a join-complete lattice embedding from
$\seq{L,\leq}$ into $\seq{\Fis L,\leq}$.
\end{lem}

The following result is much less obvious.

\begin{prop}\label{P:FisLLC}
The lattice $\seq{\Fis L,\leq}$ is lower continuous.
\end{prop}

\begin{proof}
We prove that the dual lattice $\seq{\Fis L,\subseteq}$ is upper continuous.
We put $\JJd_{i\in I}X_i=\Fls\bigl(\bigcup_{i\in I}X_i\bigr)$, for any
family $(X_i)_{i\in I}$ of elements of $\Fis L$, so it suffices to prove the
containment
  \begin{equation}\label{Eq:UppCont}
  A\cap\JJd_{i\in I}B_i\subseteq\JJd_{i\in I}(A\cap B_i),
  \end{equation}
for any $A\in\Fis L$, any upper directed poset $I$, and any
increasing (for the inclusion) family $(B_i)_{i\in I}$ of elements of $\Fis L$.
Put $B=\bigcup_{i\in I}B_i$, observe that $B$ is a filter of $L$; thus, by
Lemma~\ref{L:MonClass}, $\Cls(B)=\JJd_{i\in I}B_i$.
Denote by $C$ the right hand side of \eqref{Eq:UppCont}, and put
  \[
  D=\bigcap_{a\in A}(\JI{C}{a}).
  \]
For all $i\in I$ and $b\in B_i$, the element $a\vee b$ belongs to $A\cap B_i$,
thus to $C$; whence $B\subseteq D$.
It follows from Lemma~\ref{L:avnX} that $D$
is $\BS$-closed, thus $\Cls(B)\subseteq D$.
This means that $a\vee b\in C$, for
all $a\in A$ and all $b\in\Cls(B)=\JJd_{i\in I}B_i$, which concludes the proof
of~\eqref{Eq:UppCont}.
\end{proof}

\section{Join-semidistributive lattices of filters}
\label{S:SD+Fil}

The main result of the present section invokes the \emph{special notions of
convergence} introduced in Definition~\ref{D:CovNot}.

\begin{prop}\label{P:SD+Fil}
Let $\BS$ be a special notion of convergence on a lattice $L$. Then the
lattice $\Fis L$ \pup{with reverse inclusion} is complete, lower
continuous, and \jsd.
\end{prop}

\begin{proof}
We have already observed that $\Fis L$ is complete, and, by
Proposition~\ref{P:FisLLC}, lower continuous. In order to prove that $\Fis L$
is \jsd, it suffices to prove that for all $A$, $B$, $C$, $D\in\Fis L$ such
that $A\cap B=A\cap C=D$, the containment $A\cap(B\jjd C)\subseteq D$ holds,
where we put $B\jjd C=\Fls(B\cup C)$. Let $F$ denote the filter of $L$
generated by $B\cup C$.

We prove that $a\vee x\in D$, for all $\seq{a,x}\in A\times F$. Let
$\seq{b,c}\in B\times C$ such that $b\wedge c\leq x$. We prove, by
induction on $|\bs|$,
that $a*\bs\in D$, for all $\bs\in\fin{\set{b,c}}\setminus\set{\es}$. For
$\bs=\seq{b}$, we have $a*\bs=a\vee b\in D$, and similarly for $\bs=\seq{c}$.
Suppose that $|\bs|>1$. It follows from the induction hypothesis that
$a*\bs_*\in D$, thus $e(s)\wedge(a*\bs_*)\in B\cup C$, and therefore
  \[
  a*\bs=a\vee(e(s)\wedge(a*\bs_*))\in D.
  \]
Since the set $D$ is $\BS$-closed and $\BS$ is a special notion of
convergence, 
$a\vee(b\wedge c)=\bigwedge\famm{a*\bs}{\bs\in\fin{\set{b,c}}}$
belongs to $D$. Since $D$ is an upper subset of $L$, it follows that
$a\vee x\in D$. Hence we have proved the containment $F\subseteq E$, where
we put
  \[
  E=\bigcap\famm{\JI{D}{a}}{a\in A}.
  \]
Since $E$ is $\BS$-closed, we conclude, by using Lemma~\ref{L:MonClass}, that
$B\jjd C=\Fls(F)=\Cls(F)\subseteq E$. This means
that $a\vee x\in D$, for all $a\in A$ and all $x\in B\jjd C$.
\end{proof}

We immediately obtain the following corollary, see
\cite[Proposition~2.18]{Erne88}. Compare with Example~\ref{Ex:FilFL(3)}.

\begin{cor}\label{C:LocFinSD+}
The filter lattice $\Fil L$ of every locally finite, \jsd\ lattice $L$ is \jsd.
\end{cor}

\begin{proof}
Let $\BS$ denote the set of all finite subsets of $L$ with a least element. It
is obvious that $\BS$ is a special notion of convergence on $L$, furthermore
$\Fis L=\Fil L$ is a dually algebraic lattice. The conclusion follows from
Proposition~\ref{P:SD+Fil}.
\end{proof}

We recall that there exists a \jsd\ lattice that cannot be embedded into any
complete, \jsd\ lattice, see \cite[Example~3.25]{AGT}. Whether or not such a
lattice can be taken finitely generated is apparently an open problem, see
\cite[Problem~2]{AGT}.

\section{Sublattices of complete, lower continuous, \jsd\
lattices}\label{S:SubCLCSD+}

In this section we shall reap the consequences of
Sections~\ref{S:AxSDI}--\ref{S:SD+Fil}. We first prove a simple lemma.

\begin{lem}\label{L:CompVar}
Let $\kappa$ be either an infinite cardinal or $\infty$.
Let $K$ be a sublattice of a $\kappa$-lower continuous lattice~$L$. Denote by
$K'$ the set of all meets of downward directed subsets of $K$ with at most
$\kappa$ elements. Then~$K'$ is a sublattice of $L$, and it belongs to the
same variety as $K$.
\end{lem}

\begin{proof}
The fact that $K'$ is a sublattice of $L$ follows from the continuity
assumption on $L$. Now let $m>0$ and let $\xs$ and $\xt$ be lattice terms with
$m$ variables such that~$K$ satisfies the identity $\xs=\xt$. Let $a_0$,
\dots, $a_{m-1}\in K'$. For all $i<m$, there exists a downward directed
subset $X_i$ of $K$ with at most $\kappa$ elements such that
$a_i=\bigwedge X_i$. By using the continuity assumption on $L$, we obtain:
  \begin{align*}
  \xs(a_0,\dots,a_{m-1})&=\bigwedge\bigl(\xs(x_0,\dots,x_{m-1})\mid
  x_i\in X_i,\text{ for all }i<m\bigr)\\
  &=\bigwedge\bigl(\xt(x_0,\dots,x_{m-1})\mid
  x_i\in X_i,\text{ for all }i<m\bigr)\\
  &=\xt(a_0,\dots,a_{m-1}),
  \end{align*}
which proves that $K'$ satisfies the identity $\xs=\xt$.
\end{proof}

\begin{thm}\label{T:SubCLCSD+}
For any lattice $L$, the following statements are equivalent:
\begin{enumerate}
\item $L$ has a join-complete lattice embedding into some
complete, lower continuous, \jsd\ lattice that belongs to the same variety
as~$L$.

\item $L$ has a lattice embedding into some complete,
lower continuous, \jsd\ lattice.

\item $L$ has a lattice embedding into some $\aleph_0$-meet-complete,
$\aleph_0$-lower continuous, \jsd\ lattice.

\item $L$ satisfies the axiom \SDI.
\end{enumerate}
\end{thm}

\begin{proof}
(i)$\Rightarrow$(ii) and (ii)$\Rightarrow$(iii) are trivial.

(iii)$\Rightarrow$(iv) follows immediately from
Propositions~\ref{P:StD}, \ref{P:SD+2SDI}, and \ref{P:PresSDI}(i).

(iv)$\Rightarrow$(i)
Let $L$ be a lattice satisfying \SDI. We put
  \[
  U(\bs;a,b,c)=\setm{\bs\loz(a*\bt)}
  {\bt\in\fin{\set{b,c}}\setminus\set{\es}},
  \]
for all $a$, $b$, $c\in L$ and all $\bs\in\fin L$. Furthermore, we put
  \[
  \BS=\setm{U(\bs;a,b,c)}{a,\,b,\,c\in L\text{ and }\bs\in\fin L}.
  \]
It follows from the assumption \SDI\ that $\BS$ satisfies (S1). Let
$X\in\BS$ and $u\in L$. Write $X=U(\bs;a,b,c)$, for some $a$,
$b$, $c\in L$ and $\bs\in\fin L$. It follows from Lemma~\ref{L:circjm}
that $u\vee X=U(j(u,\bs);a,b,c)$ belongs to $\BS$, whence $u\vee X$ has a
meet in~$L$, and consequently, by using \SDI\ and Lemma~\ref{L:circjm},
  \[
  \bigwedge(u\vee X)=j(u,\bs)\loz(a\vee(b\wedge c))
 =u\vee(\bs\loz(a\vee(b\wedge c)))=u\vee\bigwedge X.
  \]
On the other hand, observe that $b\wedge c\leq a*\bt$, for all
$\bt\in\fin{\set{b,c}}$, whence
$u\wedge X=U(m(u,b\wedge c,\bs);a,b,c)$ belongs to $\BS$. Therefore, $\BS$
satisfies (S2) and (S3). Finally,
$\setm{a*\bt}{\bt\in\fin{\set{b,c}}\setminus\set{\es}}$ is equal to
$U(\es;a,b,c)$, thus it belongs to~$\BS$ and its meet is $a\vee(b\wedge c)$,
and hence $\BS$ is a special notion of convergence on $L$.

By Lemma~\ref{L:CanEmbS}, $L$ has a join-complete lattice embedding into
$\hL=\Fis L$. By the paragraph above and by Proposition~\ref{P:SD+Fil}, $\hL$
is complete, lower continuous, and \jsd. For a subset $X$ of $\hL$, let
$X^{\dnw\wedge}$ denote the set of all meets of downwards directed subsets
of~$X$. Put $L_0=L$, $L_{\xi+1}=(L_{\xi})^{\dnw\wedge}$ for any ordinal
$\xi$, and $L_{\lambda}=\bigcup_{\xi<\lambda}L_{\xi}$ for any limit ordinal
$\lambda$. It follows from Lemma~\ref{L:CompVar} that $L_{\xi}$ is a
sublattice of $\hL$ (thus it is \jsd), and it belongs to the same
variety as $L$, for any ordinal $\xi$. Hence the
same holds for $L^*=\bigcup_{\xi}L_{\xi}$. Furthermore, $L^*$ is a complete
meet-subsemilattice of $\hL$. Therefore, $L^*$ is complete, lower
continuous, \jsd, and the inclusion map from~$L$ into $L^*$ is
join-complete.
\end{proof}

\begin{cor}\label{C:SDI2DSD}
Every lattice satisfying \SDI\ is dually staircase $*$-distributive.
\end{cor}

\begin{proof}
Let $L$ be a lattice satisfying \SDI. It follows from Theorem~\ref{T:SubCLCSD+}
that $L$ can be embedded into some complete, lower continuous, and \jsd\
lattice $\ol{L}$. It follows from Proposition~\ref{P:StD} that $\ol{L}$ is
dually staircase distributive. Hence, by Lemma~\ref{L:SubL*}, $L$ is
also dually
staircase distributive.
\end{proof}

Hence, by using Proposition~\ref{P:SD+2SDI}, we obtain the following.

\begin{cor}\label{C:ZiStD}
Let $L$ be a $\aleph_0$-lower continuous lattice. If $L$ is dually zipper
distributive, then $L$ is dually staircase distributive.
\end{cor}

We emphasize that we do \emph{not} assume completeness of $L$ in the statements
of Corollaries~\ref{C:SDI2DSD} and~\ref{C:ZiStD}.

By using Corollary~\ref{C:DNoethSDI}, we obtain the following.

\begin{cor}\label{C:LBandSD+lf}
Every lower bounded lattice can be embedded into some complete, lower
continuous, \jsd\ lattice.
\end{cor}

One cannot hope to strengthen the conclusion of Corollary~\ref{C:LBandSD+lf} by
requiring the larger lattice to be lower bounded; see also
\cite[p.~207]{AdGo01}. The proof of the following lemma is straightforward,
and left to the reader.

\begin{lem}\label{L:SubLB}
Let $K$ be a sublattice of a lattice $L$. We assume that the inclusion map
from $K$ into $L$ is lower bounded. If a nonempty subset $X$ of $K$
has a meet in~$L$, then $X$ has a meet in $K$, and the two meets are equal.
\end{lem}

\begin{exple}\label{Ex:NonCplLB}
The free lattice on three generators $\FL(3)$ cannot be embedded into any
$\aleph_0$-meet-complete, lower bounded lattice.
\end{exple}

\begin{proof}
Suppose that $\FL(3)$ is a sublattice of a $\aleph_0$-meet-complete, lower
bounded lattice $L$. Since $\FL(3)$ is finitely generated and $L$ is lower
bounded, the inclusion map $f\colon\FL(3)\hookrightarrow L$ is lower bounded.
Hence, it follows from Lemma~\ref{L:SubLB} that $\FL(3)$ is
$\aleph_0$-meet-complete, which is known not to be the case (see
\cite[Section~I.5]{FJN}).
\end{proof}

It is noteworthy to record the following immediate consequence of
Proposition~\ref{P:PresSDI}, Corollary~\ref{C:FinGenSDI}, and
Theorem~\ref{T:SubCLCSD+}.

\begin{cor}\label{C:FinSDIEquiv}
A lattice $L$ has an embedding into some complete, lower continuous, \jsd\
lattice if{f} every finitely generated sublattice of $L$ has such an embedding.
\end{cor}

The following example is the lower part of Example~3.25 in \cite{AGT}, and it
is also the dual, minus the top element, of the lattice $Z_\omega$ of
\cite[Page~299]{Erne88}. It shows that lower continuity cannot be replaced by
upper continuity in the statement of Corollary~\ref{C:LBandSD+lf}.

\begin{exple}\label{Ex:NonEmbLB}
The lattice $L$ of Figure~\textup{1} is locally finite and lower bounded, but
it cannot be embedded into any $\aleph_0$-join-complete,
$\aleph_0$-upper continuous, \jsd\ lattice.
\end{exple}

\begin{figure}[htb]
\includegraphics{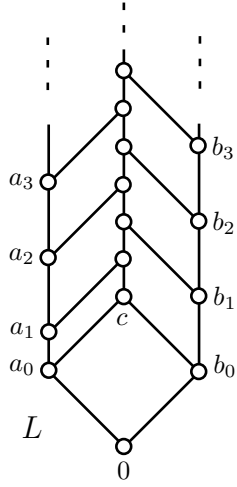}
\caption{A locally finite lower bounded lattice.}
\end{figure}

\begin{proof}
It is straightforward to verify that the lattice $L$ is locally finite
and lower bounded. Let $L'$ be a $\aleph_0$-join-complete, $\aleph_0$-upper
continuous, \jsd\ lattice containing $L$. Put $a=\bigvee_{n<\omega}a_n$ and
$b=\bigvee_{n<\omega}b_n$. Since $a_n\leq b_n\vee c$ and
$b_n\leq a_{n+1}\vee c$, for all $n<\omega$, the equality $a\vee c=b\vee c$
holds. Since $L'$ is \jsd, the inequality
$a\leq(a\wedge b)\vee c$ holds, thus, \emph{a fortiori},
$a_1\leq(a\wedge b)\vee c$. However, since $L'$ is $\aleph_0$-join complete and
$\aleph_0$-upper continuous, $a\wedge b=\bigvee_{n<\omega}(a_n\wedge
b_n)=0$, so
we obtain the inequality $a_1\leq c$, \contr.
\end{proof}

On the other hand, the duals of \cite[Proposition~2.18]{Erne88} and
\cite[Lemma~10]{Rein95} imply immediately the following example. It illustrates
the importance of the choice of $\BS$ in the proof of
Corollary~\ref{C:LBandSD+lf}. Compare also with Corollaries~\ref{C:Id(LB)} and
\ref{C:LocFinSD+}.

\begin{exple}\label{Ex:FilFL(3)}
The filter lattice $L$ of $\FL(3)$ is not \jsd.
\end{exple}

We conclude this section by the following analogue of
Theorem~\ref{T:SubCLCSD+} for \emph{complete} embeddings.

\begin{figure}[htb]
\includegraphics{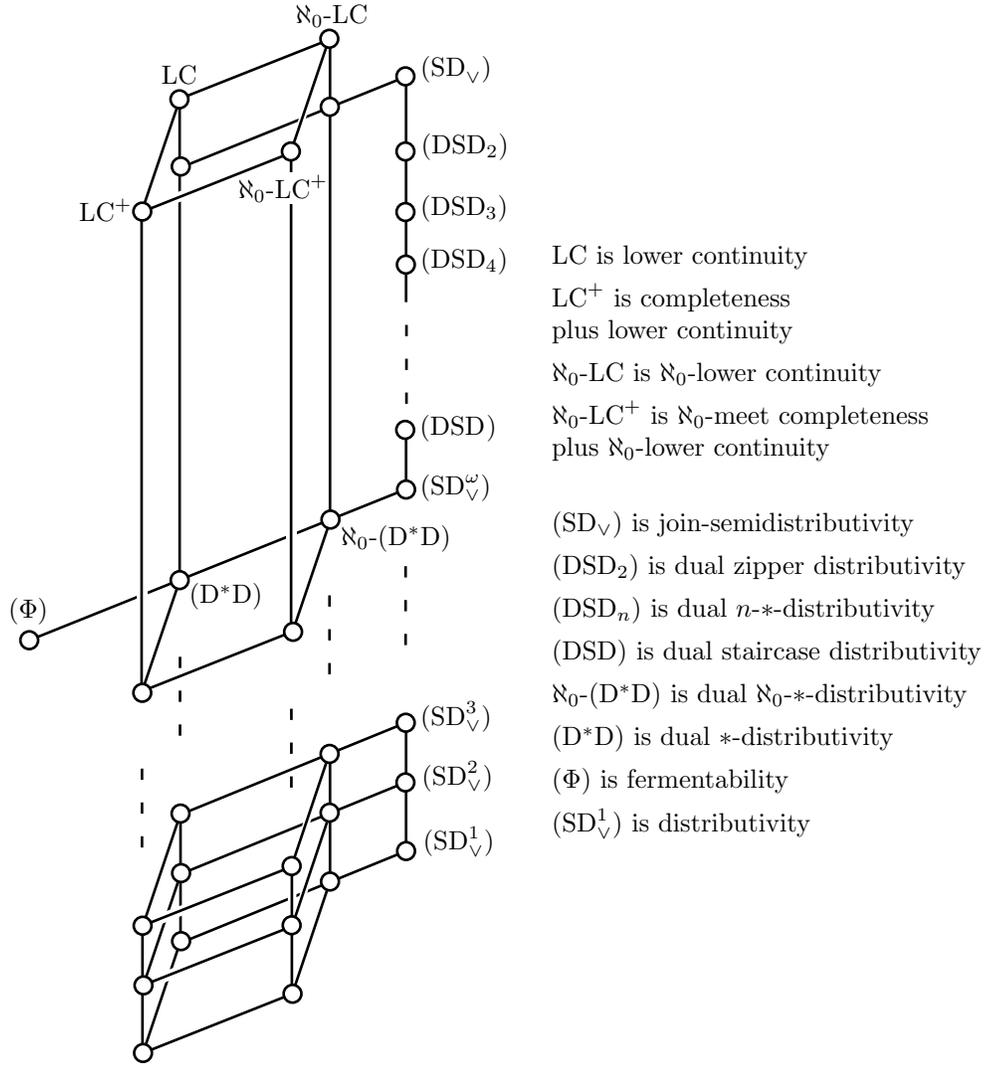}
\caption{Classes of lattices.}
\end{figure}

\begin{thm}\label{T:CpSubCLCSD+}
For any lattice $L$, the following statements are equivalent:
\begin{enumerate}
\item $L$ has a complete lattice embedding into some complete, lower
continuous, \jsd\ lattice.

\item $L$ has a meet-complete lattice embedding into some complete, lower
continuous, \jsd\ lattice.

\item $L$ is dually $*$-distributive.

\item $L$ is dually zipper distributive and lower continuous.
\end{enumerate}
\end{thm}

\begin{proof}
(i)$\Rightarrow$(ii) is trivial.

(ii)$\Rightarrow$(iii) follows immediately from Proposition~\ref{P:StD} and
Lemma~\ref{L:SubL*}.

(iii)$\Rightarrow$(iv) follows immediately from Proposition~\ref{P:StD}.

(iv)$\Rightarrow$(i). Let $L$ be dually zipper distributive and lower
continuous. We put
  \[
  \BS=\setm{X\subseteq L}{X\text{ is downward directed and }
  \bigwedge X\text{ exists}}.
  \]
It is straightforward to verify that $\BS$ is a special notion of convergence
on $L$ and that the map $x\mapsto\upw x$ is a complete lattice embedding from
$L$ into $\Fis L$. By Proposition~\ref{P:SD+Fil}, $\Fis L$ is complete, lower
continuous, and \jsd.
\end{proof}

Most classes of lattices encountered in the present paper are
represented on Figure~2. The largest classes are on the top of the diagram. The
diagram without $(\Phi)$ (i.e., fermentability) is a meet-semilattice, for
example, the intersection of $(\mathrm{SD}_\vee)$ and $(\mathrm{LC}^+)$ is,
indeed, contained in the class $(\mathrm{D}^*\mathrm{D})$ of all dually
$*$-distributive lattices.

\section{Non-embeddability results into bi-algebraic lattices}\label{S:NonEmb}

We first state the central lemma underlying all the results of the present
section.

\begin{lem}\label{L:NonHom}
Let $L$ be a $\aleph_0$-complete, $\aleph_0$-upper continuous, and
$\aleph_0$-lower continuous lattice with zero. Let $(a_n)_{n<\omega}$ be a
sequence of elements of $L$ and let $c\in L$ such that the following statements
hold:

\begin{enumerate}
\item
$\left(\bigvee_{i<m}a_i\right)\wedge\left(\bigvee_{j<\omega}a_{m+j}\right)=0$,
for all $m<\omega$.

\item $a_0\leq\bigvee_{j<\omega}a_{m+j}\vee c$, for all $m<\omega$.

\item $a_0\wedge c=0$.
\end{enumerate}

Then $a_0=0$.
\end{lem}

\begin{proof}
Put $b_n=\bigvee_{j<\omega}a_{n+j}$, for all $n<\omega$, and
$b=\bigwedge_{n<\omega}b_n$. It follows from (i) that
$\left(\bigvee_{i<n}a_i\right)\wedge b=0$, for all $n<\omega$, hence, by the
$\aleph_0$-upper continuity of $L$, we obtain that $b=0$. Since (ii) can be
written $a_0\leq b_m\vee c$ for all
$m<\omega$, it follows from the $\aleph_0$-lower continuity of $L$ that
$a_0\leq b\vee c=c$. Therefore, from (iii) it follows that $a_0=0$.
\end{proof}

\begin{cor}\label{C:NonHom}
Let $L$ be a lattice with zero, let $(a_n)_{n<\omega}$ be a
sequence of elements of $L$, and let $c\in L$ such that $a_0\neq0$ and
the following statements hold:
\begin{enumerate}
\item $\left(\bigvee_{i<m}a_i\right)\wedge\left(\bigvee_{j<n}a_{m+j}\right)=0$,
for all $m$, $n<\omega$.

\item For all $m<\omega$, there exists $n<\omega$ such that
$a_0\leq\bigvee_{j<n}a_{m+j}\vee c$.

\item $a_0\wedge c=0$.
\end{enumerate}

Then $L$ cannot be embedded into any $\aleph_0$-complete, $\aleph_0$-upper
continuous, and $\aleph_0$-lower continuous lattice.
\end{cor}

\begin{proof}
Let $L'$ be a $\aleph_0$-complete, $\aleph_0$-upper continuous, and
$\aleph_0$-lower continuous lattice such that $L$ embeds into $L'$.
After replacing $L'$ by $\upw L$, we
may assume, without loss of generality, that $L'$ has the same zero as $L$.
It is trivial that (ii) and (iii) of the statement of Corollary~\ref{C:NonHom}
imply, respectively, (ii) and (iii) of the statement of Lemma~\ref{L:NonHom}.
Since $L'$ is $\aleph_0$-upper continuous, (i) follows as well. By
Lemma~\ref{L:NonHom}, $a_0=0$, a contradiction.
\end{proof}

As usual, in any modular lattice $L$ with zero, let $c=a\oplus b$ mean that
$c=a\vee b$ while $a\wedge b=0$. A family $(a_i)_{i\in I}$ of elements of
$L$ is \emph{independent}, if the equality
  \[
  \left(\bigvee\famm{a_i}{i\in X}\right)\wedge
  \left(\bigvee\famm{a_i}{i\in Y}\right)=
  \bigvee\famm{a_i}{i\in X\cap Y}
  \]
holds, for all finite subsets $X$ and $Y$ of $I$. Since $L$ is
modular, it is sufficient to verify this for $X$ a
singleton, see \cite{Grat98}. We say that $a$, $b\in L$ are
\emph{perspective}, if there exists $c\in L$ such that $a\oplus c=b\oplus c$.

\begin{cor}\label{C:ModHom}
Let $L$ be a modular lattice with zero, suppose that $L$ has an infinite
independent sequence of nonzero pairwise perspective elements. Then $L$ cannot
be embedded into any $\aleph_0$-complete, $\aleph_0$-upper continuous, and
$\aleph_0$-lower continuous lattice.
\end{cor}

\begin{proof}
Suppose that $L$ embeds into a $\aleph_0$-complete, $\aleph_0$-upper
continuous, and $\aleph_0$-lower continuous lattice $L'$. As above, we may
assume that $L$ and~$L'$ have the same zero. Now we use
Lemma~\ref{L:CompVar} (and its dual) as in the proof of (iv)$\Rightarrow$(i)
of Theorem~\ref{T:SubCLCSD+}, but this time by alternating $\omega_1$ times
closure under countable meets and joins. We obtain that the closure $L^*$
of~$L$ (within $L'$) under countable meets and countable joins belongs to the
same variety as $L$. In particular, in addition to being $\aleph_0$-complete,
$\aleph_0$-upper continuous, and $\aleph_0$-lower continuous, the lattice
$L^*$ is modular.

Let $(a_n)_{n<\omega}$ be a sequence of pairwise perspective elements of $L$,
with $a_0\neq0$. For all $n>0$, there exists $c_n\in L$ such that
$a_0\oplus c_n=a_n\oplus c_n=a_0\vee a_n$. It is an easy exercise to verify
that the $a_n$-s and $c=\bigvee_{n>0}c_n$
satisfy the assumptions of Corollary~\ref{C:NonHom}
(with $n=1$ in (ii)), see also the proof of \cite[Theorem~I.3.8]{Neum60} or
\cite[Satz~IV.2.1]{FMae58}. Hence, by Corollary~\ref{C:NonHom}, $a_0=0$, a
contradiction.
\end{proof}

\begin{cor}\label{C:InfVect}
Let $V$ be an infinite-dimensional vector space over a division ring. Then the
subspace lattice of $V$ cannot be embedded into any $\aleph_0$-complete,
$\aleph_0$-upper continuous, and $\aleph_0$-lower continuous lattice.
\end{cor}

\begin{cor}\label{C:Co(omega)}
Let $\seq{I,\utr}$ be an infinite chain. Then the lattice $\Co(I)$ of
order-convex subsets of $I$ cannot be embedded into any $\aleph_0$-complete,
$\aleph_0$-upper continuous, and $\aleph_0$-lower continuous lattice.
\end{cor}

\begin{proof}
Without loss of generality, $I$ has an infinite, strictly increasing sequence
$z\tr x_0\tr x_1\tr x_2\tr\cdots$. Put $a_m=\set{x_m}$, for all
$m<\omega$, and $c=\set{z}$. It is obvious that the $a_m$-s and $c$ satisfy the
assumptions (i), (ii) (for $n=1$), and (iii) of Corollary~\ref{C:NonHom}. The
conclusion follows from Corollary~\ref{C:NonHom}.
\end{proof}

Both Corollaries~\ref{C:InfVect} and \ref{C:Co(omega)} solve negatively
Problem~5 of \cite{SeWe3}.

\section{Open problems}\label{S:Pbs}

Figure~2 shows containments between various classes of \jsd\ lattices. This
suggests the following general problem.

\begin{pb}\label{Pb:IdClass}
Do all lines on Figure~2 represent proper containments?
\end{pb}

Of course, some partial answers to Problem~\ref{Pb:IdClass} are known. For
example, the classes $(\mathrm{SD}_{\vee}^n)$ and $(\mathrm{SD}_{\vee}^{n+1})$
are, indeed, distinct, for every positive integer $n$. Many related
examples can
also be found in \cite{Erne88}. On the other hand, we do
not have any example to show that dual staircase distributivity and dual zipper
distributivity are really distinct notions---they coincide in the presence of
$\aleph_0$-lower continuity, see Corollary~\ref{C:ZiStD} and Figure~2. Also
observe that a lattice $L$ is \jsd\ if{f} every three-generated sublattice of
$L$ is \jsd. A related problem, inspired by Corollary~\ref{C:FinSDIEquiv}, is
the following.

\begin{pb}\label{Pb:TestSDI}
Does there exist a positive integer $n$ such that a lattice $L$ satisfies \SDI\
if{f} every $n$-generated sublattice of $L$ satisfies \SDI?
\end{pb}

It follows from Corollary~\ref{C:LBandSD+lf} that every lower bounded lattice
can be embedded into some complete, lower continuous, \jsd\ lattice. The result
would look better if we could replace ``lower continuous'' by ``dually
algebraic'', but this we do not know.

\begin{pb}\label{Pb:LC2DAlg}
Can every lower bounded lattice be embedded into some dually algebraic \jsd\
lattice?
\end{pb}

Trying to improve the universal theory instead of the completeness
condition yields, for example, the following problem.

\begin{pb}\label{Pb:LB2LB}
Can every lower bounded lattice be embedded into some complete, lower
continuous lattice in $\mathbf{Q}(\mathcal{LB}_{\mathrm{fin}})$?
\end{pb}

It is conceivable that the extension $\Fis L$ defined in the proof of
Corollary~\ref{C:LBandSD+lf} belongs to
$\mathbf{Q}(\mathcal{LB}_{\mathrm{fin}})$, but we do not know how to
prove this.

It follows from Whitman's Theorem that every lattice can be embedded into an
algebraic and spatial lattice, namely, a partition lattice. It is
also proved in
\cite{HPR} that every \emph{modular} lattice can be embedded, \emph{within its
variety}, into an algebraic and spatial lattice.

\begin{pb}\label{Pb:DuCoOm}
Can every lattice be embedded, within its variety, into some algebraic and
spatial lattice?
\end{pb}

For every algebraic lattice $A$, the lattice $\SP(A)$ of all
algebraic subsets of $A$ (see \cite{Gorb94,Gorb}) is dually
algebraic and \jsd.

\begin{pb}\label{Pb:CplLCSPA}
Can every complete, lower continuous (resp., dually algebraic), \jsd\ lattice
be embedded into $\SP(A)$, for some complete, upper continuous (resp.,
algebraic) lattice $A$?
\end{pb}

The deepest result of \cite{AGT} is probably that every \emph{finite} \jsd\
lattice can be embedded into $\SP(A)$, for some algebraic
lattice $A$.

\section*{Acknowledgment}
I thank the anonymous referee for his careful reading of the paper, which lead
to the correction of several embarrassing oversights.

\end{document}